\documentclass[12pt]{article}

\usepackage{a4wide}
\usepackage{amssymb}
\usepackage{amsfonts}
\usepackage{amsmath}
\input xy
\xyoption{arrow} \xyoption{matrix}

\date{}

\newtheorem{proposition}{Proposition}[section]
\newtheorem{theorem}[proposition]{Theorem}
\newtheorem{lemma}[proposition]{Lemma}

\newtheorem{corollary}[proposition]{Corollary}

\def\Kdim{{\rm K.dim }\,}

\def\der{\partial }

\def\nFM0{{\nu }_{F,M_0}}
\def\nFN0{{\nu }_{F,N_0}}
\def\nGN0{{\nu }_{G,N_0}}

\def\N0{ {\bf N}_0 }

\def\t{\otimes}
\def\g{\gamma}
\def\v{\varphi}
\def\ra{\rightarrow}

\def\lra{\leftrightarrow}
\def\Xpm{X^{\pm }}

\def\s{\sigma}

\def\l1{{\lambda}_1}

\def\a{\alpha}
\def\a0{ {\alpha }_0}
\def\a1{ {\alpha }_1}

\def\l{\lambda}
\def\o{\omega}

\def\nFGM0{{\nu }_{F,G,M_0}}

\def\nFN0{{\nu}_{F,N_0}}

\def\sn{{\sigma}^n}
\def\sm{{\sigma}^m}
\def\si{{\sigma}^i}

\def\sm1{{\sigma}^{-1}}

\def\smtp1{{\sigma}^{-t+1}}

\def\o{\omega }
\def\S1{S^{-1}}

\def\Xpm1{X^{\pm 1}_1}

\def\sPM1{{\sigma }^{\pm 1}}
\def\sMP1{{\sigma }^{\mp 1 }}

\def\d{\delta}

\def\q{{\bf q}}

\def\G{\Gamma}

\def\OO{{\cal O}}

\def\Ytm1{Y^{t-1}}
\def\Yim1{Y^{i-1}}

\def\CK{{\cal K}}

\def\CF{{\cal F}}

\def\Zn{\mathbb{Z}_n}

\def\i{{\bf i}}

\def\Zn{\mathbb{Z}_n}

\def\Aut{{\rm Aut}}

\def\ad{{\rm ad }}
\def\dim{{\rm dim }}

\def\ker{ {\rm ker } }

\def\D{ \Delta }

\def\SL2Z{ {\rm SL}_2({\bf Z}) }

\def\Gp1{ G^{1 , 1 } }
\def\P11{ P^{-1 , 1 } }
\def\Pp1{ P^{1 , 1 } }

\def\Supp{{\rm Supp}}

\def\nCLsr{{}^\nu\kern-2pt {\cal L}^{\sigma , \rho  }}
\def\nP{{}^\nu \kern-2pt P}
\def\nL{{}^\nu\kern-2pt L}
\def\nLL{{}^\nu\kern-2pt \Lambda}
\def\nPsr{{}^\nu\kern-2pt P^{\sigma , \rho  }}
\def\nLsr{{}^\nu\kern-2pt L^{\sigma , \rho  }}
\def\nuCL{{}^\nu\kern-2pt  {\cal L}}
\def\nCLsr{{}^\nu\kern-2pt {\cal L}^{\sigma , \rho  }}
\def\nCL1m{{}^\nu\kern-2pt {\cal L}^{-1 , 1  }}
\def\x1nu{x^\frac{1}{\nu}}
\def\xm1nu{x^{-\frac{1}{\nu}}}

\def\bX{{\bf X}}

\def\ra{\rightarrow }

\def\coker{{\rm coker}}

\def\nAM0{{\nu }_{{\cal A},M_0}}
\def\nAN0{{\nu }_{{\cal A},N_0}}

\def\Kdim{ {\rm Kdim } }

\def\ad{ {\rm ad }}

\def\gn{\mathfrak{n}}
\def\gm{\mathfrak{m}}
\def\gp{\mathfrak{p}}
\def\gq{\mathfrak{q}}

\def\j{{\bf j}}

\def\Fq{\mathbb{F}_q}
\def\gldim{{\rm gl.dim}}
\def\hC{\widehat{C}}
\def\bY{\overline{Y}}
\def\bX{\overline{X}}
\def\bC{\overline{C}}
\def\Max{{\rm Max}}
\def\ann{{\rm ann}}
\begin{document}

\author{V.\  Bavula}

\title{The Carlitz Algebras}

\maketitle
\begin{abstract}
The Carlitz $\mathbb{F}_q$-algebra $C=C_\nu$, $\nu \in
\mathbb{N}$, is generated by an algebraically closed field $\CK $
(which contains a non-discrete locally compact field of positive
characteristic $p>0$, i.e. $K\simeq \mathbb{F}_q[[ x,x^{-1}]]$,
$q=p^\nu$), by the (power of the) {\em Frobenius} map $X=X_\nu
:f\mapsto f^q$, and by the {\em Carlitz derivative} $Y=Y_\nu$. It
is proved that the
 Krull  and  global dimensions of $C$ are  $2$, a classification of
simple $C$-modules and ideals are given, there are only {\em
countably many} ideals, they commute $(IJ=JI)$, and each ideal is
a unique product of maximal ones.  It is a remarkable fact that
any simple $C$-module is a sum of eigenspaces of the element $YX$
(the set of eigenvalues for $YX$ is given explicitly for each
simple $C$-module). This fact is crucial in finding the group
$\Aut_{\Fq}(C)$ of $\Fq$-algebra automorphisms of $C$ and in
proving that  two distinct Carlitz rings are not isomorphic
$(C_\nu \not\simeq C_\mu$ if $\nu \neq \mu$). The centre of $C$ is
found explicitly, it is a UFD that contains {\em countably many}
elements.

\end{abstract}


\section{Introduction}
In this paper, module means a {\em left} module. Recall that
$A_1(F)=F\langle X, \der \, | \, \der X-X\der =1\rangle$ is the
{\em first Weyl} algebra over a field $F$.  Let $k$ and $l=l_p$ be
fields of characteristic zero and $p>0$ respectively. The first
Weyl algebras $A_1(k)$ and $A(l)$ have very distinctive
properties, name just a few: $A_1(k)$ is a simple algebra but
$A(l)$ is not (there are finite and infinite dimensional factor
algebras of $A(l)$); $A_1(k)$ has a small centre (equal to $k$)
but $A(l)$ has a big centre (equal to a polynomial algebra $l[X^p,
\der^p]$ in two variables); the Krull and global dimensions $\Kdim
(A_1(k))=\gldim (A_1(k))=1$ but $\Kdim (A_1(l))=\gldim
(A_1(l))=2$; the classical Krull dimension ${\rm cl.Kdim}
(A_1(k))=0$ but ${\rm cl.Kdim}(A_1(l))=2$; the algebra $A_1(k)$
has no nonzero simple finite dimensional modules but all the
simple modules over the algebra  $A(l)$ are finite dimensional;
only  a tiny bit of simple $A_1(k)$-modules are weight modules but
all simple $A_1(l)$-modules are weight and finite dimensional over
$l$.

The Carlitz algebras $C_\nu $ exhibit properties that are ``in
between'' of that of $A_1(k)$ and $A_1(l)$: the centre of $C_\nu$
is very small (it contains countably many elements), $\Kdim (C_\nu
)=\gldim (C_\nu )=2$ and
 ${\rm cl.Kdim} (C_\nu)=1$, all simple $C_\nu$-modules are weight
 but there are plenty of simple infinite dimensional and finite
 dimensional $C_\nu$-modules over $\CK $, every factor ring of $C$
 is infinite dimensional over $\CK$.

The paper has the following structure: In Section \ref{C=GWA}, we
recall the definition of the Carlitz algebras and prove that they
are generalized Weyl algebras. Using this fact, in Section
\ref{simC}, a classification of simple $C$-modules is given using
some results of \cite{BavGWA92} and \cite{BavVOYJA97} (Theorems
\ref{wlins} and \ref{Cscyc}). The Carlitz algebras have pleasant
representation theory due to the fact that every simple module is
{\em weight} (i.e. a direct sum of eigenspaces for the element
$YX$, Corollary \ref{1crfn}). Any nonunit element of $C$ has {\em
finite dimensional} kernel and cokernel over $\CK$ on any simple
$C$-module (Theorem \ref{fdkck}).

In Section \ref{prCarl}, the ideals of $C$ are explicitly
described (Theorem \ref{mCarid}), name just a few (rather
peculiar) results: ideals commute $(IJ=JI)$, and there are
countably many of them, every ideal is a unique product of maximal
ones (as in the case of the ring of integers), each simple factor
ring is {\em infinite dimensional} over $\CK$ and has Krull and
global dimension $1$. The centre of $C$ is found (Lemma
\ref{Rs2}.(3)), it is a UFD and contains countably many elements.

In Section \ref{1kglC}, it is proved that two distinct Carlitz
algebras are not isomorphic (Theorem \ref{caris}) and the groups
$\Aut (C)$ and $\Aut_\CK (C)$ of ring automorphisms and of
$\CK$-automorphisms are found explicitly (Theorem \ref{grrautC}).

In Section \ref{kglC}, it is proved that the ring $C$ has Krull
and global dimension $2$.


\section{The Carlitz algebras are generalized Weyl algebras}\label{C=GWA}

{\bf Generalized Weyl algebras}.  Let $D$ be a ring with an
automorphism $\s $ and a central element $a$.
 A {\bf generalized Weyl algebra} (a GWA, for short) $A=D(\sigma, a)$  of degree 1,
is the ring generated by $D$ and two indeterminates $X$ an $Y$
subject to the relations \cite{BavFA91}:
$$
X\alpha=\sigma(\alpha)X \ {\rm and}\
Y\alpha=\sigma^{-1}(\alpha)Y,\; {\rm for \; all }\; \alpha \in D,
\ YX=a \ {\rm and}\ XY=\sigma(a).
 $$
The algebra $A={\oplus}_{n\in \mathbb{Z}}\, A_n$ is a
$\mathbb{Z}$-graded algebra where $A_n=Dv_n$,
 $v_n=X^n\,\, (n>0), \,\,v_n=Y^{-n}\,\, (n<0), \,\,v_0=1.$
 It follows from the above relations that
$$v_nv_m=(n,m)v_{n+m}=v_{n+m}<n,m> $$
for some $(n,m)=\s^{-n-m}(<n,m>)\in D$. If $n>0$ and $m>0$ then
$$n\geq m:\, (n,-m)=\sigma^n(a)\cdots \sigma^{n-m+1}(a),\, (-n,m)=\sigma^{-n+1}(a)\cdots \sigma^{-n+m}(a),$$
$$n\leq m:\,\,\,\, (n,-m)=\sigma^{n}(a)\cdots \sigma(a),\,\,\,
(-n,m)=\sigma^{-n+1}(a)\cdots a,$$ in other cases $(n,m)=1$.

Let $K[H]$ be a polynomial ring in one variable $H$ over the field
$K$, $\s :H\ra H-1$ be the $K$-automorphism of the algebra $K[H]$
and $a=H$. The first Weyl algebra $A_1$ is isomorphic to the
generalized Weyl algebra
$$A_1\simeq K[H](\s , H),\; X\lra X,\; \der\lra Y,\; \der X\lra H.$$

{\bf The Carlitz algebras}. Any non-discrete locally compact field
of positive characteristic $p>0$ is isomorphic to the field $K=\Fq
[[x,x^{-1}]]$ of Laurent series with coefficients from the Galois
field $\Fq$ that contains $q:=p^\nu$ elements $\nu \in
\mathbb{N}$. A typical nonzero element $\l $ of $K$ is a series
$\sum_{i=m}^\infty \l_ix^i$ where $m\in \mathbb{Z}$, $\l_i\in
\Fq$, $\l_m\neq 0$. The field $K$ is equipped with a {\em
non-archimedian absolute value} $|\l  |:=q^{-m}$ which can be
extended onto the  the completion of an algebraic closure
$\overline{K}$ of the field $K$. Let $\CK$ be an algebraically
closed field extension of $K$. We denote by $F:a\mapsto F(a):=a^p$
the {\em Frobenious map} of $\CK$ ($K$ or $\Fq$), and let $\s
=\s_\nu :=F^\nu$ where $q:=p^\nu$. The fields $\CK$, $K$
 and $\Fq$ are topological fields with respect to the $p$-adic
 filtration induced by $| \cdot |$, and the maps $F$ and $\s $ are
 continuous. The action of $F$ and $\s $ on $\CK$ can be
 extended to the action on the set of maps $\CF $ from $K$ to
 $\CK$ by the rule $F:f\mapsto f^p$ and $X:f\mapsto f^q$ where
 $X=X_\nu$ denotes the extension of the $\s $ (we have introduced
 another letter in order to avoid confusion later).  There are two
 distinguished maps $\D $ and $Y$ from $\CF $ to $\CF $
 \cite{Carlitz1935}:
 \begin{eqnarray*}
{\bf the\;\;  difference \;\; operator}: &  \D u(t):= u(xt)-xu(t), \;\; u(t)\in \CF, \\
 {\bf the \;\;  Carlitz \;\; derivative}: & Y:=\sqrt[q]{\D}.
\end{eqnarray*}

The {\bf Carlitz algebra} $C=C_\nu$ is the subalgebra of ${\rm
End}_{\Fq} (\CF )$ generated by the isomorphic copy of the field
$\CK$ in the endomorphism algebra ${\rm End}_{\Fq} (\CF )$ ($\CK
\ra {\rm End}_{\Fq} (\CF )$, $ \l \mapsto (u\mapsto \l u)$) and
the maps $X=X_\nu$ and $Y=Y_\nu$.

The elements $X$ and $Y$ satisfy the following relations
\cite{Koch98}:
$$ YX-XY=[1]^{\frac{1}{q}}, \;\; X\l =\l^qX, \;\; Y\l
=\l^{\frac{1}{q}}Y\;\; (\l \in \CK ),$$ where $[1]:= -x+x^q$. The
algebra $C$ is a Noetherian domain and any element of the algebra
$C$ is a unique sum $\sum \l_{ij}X^iY^j$ where $\l_{ij}\in \CK$
\cite{KochJNT00}. It follows immediately that the map
$$ C\ra \CK [H](\s, H), \;\; X\mapsto X, \;\; Y\mapsto Y, \;\;
YX\mapsto H,\;\;  \l \mapsto \l \;\; (\l \in \CK ), $$ is an
$\Fq$-algebra isomorphism (as the GWA  $\CK [H](\s,
H)=\oplus_{i\in \mathbb{Z}}\CK [H]v_i$) where $$ \s \in \Aut (\CK
[H]): \, H\mapsto H-\l_1, \;\; \l \mapsto \l^q \; (\l \in \CK),$$
and  $\l_1:=[1]^\frac{1}{q}=-x^\frac{1}{q}+x$. So, the Carlitz
algebra is the generalized Weyl algebra $C=D(\s , H)$ with
coefficients from the polynomial algebra $D:=\CK [H]$ with
coefficients from the field $\CK$. Since the polynomial algebra
$D$ is a Noetherian domain then so is the GWA $C$, by
\cite{BavGWA92}, Proposition 1.3. Note that the Carlitz algebra is
an algebra over the {\em finite} field $\mathbb{F}_q$ and is {\em
not} an algebra over $\CK$ or even over $\mathbb{F}_{q^i}$, $i\geq
2$.


\section{A classification of simple modules over the Carlitz
algebras}\label{simC}

The representation theory of GWAs are well understood, see for
example \cite{BavGWA92} and \cite{BavVOYJA97}. We use  some of the
results of these two papers to give a classification of simple
modules over the Carlitz algebras. A surprising feature is that
all simple $C$-modules are weight modules (that is $H$ acts as a
semi-simple linear map).

Let $\CK (H)=\S1 D$ the field of fractions of $D$ where
$S=D\backslash \{ 0\}$. The localization $B=\S1 C$ of
 the ring $C$ at $S$ is the skew Laurent polynomial ring $B=\CK (H)[X, X^{-1}; \s ]$
  with coefficients from the field $\CK (H)$.
  We may identify the ring $C$ with a subring of $B$ via the ring monomorphism
  $$ C\ra B, \; X\mapsto X,\; Y\mapsto HX^{-1}, \; d\mapsto d \; (d\in D).$$
  The ring $B$ is an Euclidean ring (a left and right division algorithm with remainder holds), hence a principal left and right ideal
domain.
   We have the partition of the set $\hC$ of isoclasses of simple
   $C$-modules
   \begin{equation}\label{Aunion}
   \hC =\hC (D-{\rm torsion})\cup \hC (D-{\rm torsionfree}),
\end{equation}
where a simple $C$-module $M$ is $D$-{\em torsion} (respectively,
$D$-{\em torsionfree}) if $\S1 M=0$ (respectively, $\S1 M\ne 0$).
We will see shortly that $\hC (D- {\rm torsionfree})=\emptyset$
(Corollary \ref{1crfn}).

{\bf Max($D$).} Let $G=\langle \s \rangle$ be the subgroup of the
group of ring automorphisms $\Aut (D)$ of $D$ generated by the
element $\s $. The group $G$  acts in the  obvious way on the set
of maximal ideals of the algebra $D$,
$${\rm Max}(D):= \{ D(H-\l )\, | \, \l \in \CK\}\simeq \CK, \;\; D(H-\l )\lra \l .  $$
 So, the  orbit $\OO $ of an
element $\gp\in {\rm Max}(D)$ is   $\OO (\gp ) =\{\s {}^i(\gp
),i\in \mathbb{Z}\}$. In more detail, if $\gp=D(H-\l )$ for some
$\l \in \CK$ then an easy induction gives that for each $\l \in
\CK$ and each natural number $n\geq 1$: 
\begin{equation}\label{spnHl}
\sn (H-\l )=H-\l_1-\l_1^q-\cdots -\l_1^{q^{n-1}}-\l^{q^n},
\end{equation}
\begin{equation}\label{smnHl}
\s^{-n} (H-\l )=H+\l_1+\l_1^\frac{1}{q}+\cdots
+\l_1^{\frac{1}{q^{n-1}}}-\l^{\frac{1}{q^n}}.
\end{equation}
An orbit is called {\em cyclic} of length $n$ (respectively, {\it
linear}) if it contains a {\em finite} (respectively, {\em
infinite}) number $n=|\OO |$ of elements. The set of cyclic (resp.
linear) orbits is denoted by {\it Cyc}
 (resp. {\it Lin}). An orbit $\OO $ is  called {\em  degenerate}, if it contains
a maximal ideal $\gp $ such that $H\in \gp $ (such ideals are
called {\em marked}). We denote by {\it Cycd} and {\it Lind}
(resp. {\it Cycn} and {\it Linn}) the set of all degenerate (resp.
non-degenerate) cyclic and linear orbits, respectively.

Each linear orbit $\OO (\gp ) $, via the map $\OO (\gp ) \ra
\mathbb{Z},\s {}^i(\gp)\ra i,$ may be identified with the set of
integers $\mathbb{Z}$. Therefore, for $\OO (\gp ) \in $ {\it Lin}
we may use, without mentioning it explicitly, all the definitions
and notations which are employed for $\mathbb{Z}$ (such as the
order, the segment, the interval, etc.). For example, $\s {}^i(\gp
)\leq \s {}^j(\gp )$ iff $i\leq j;$ $(-\infty ,\gp ]:=\{\s
{}^i(\gp ),i\leq 0\}$, etc.
 Marked ideals $\gp {}_1
<\cdots <\gp {}_s$  of a degenerate linear  orbit $\OO $ divide it
into $s+1$  parts,
$$\G {}_1=(-\infty ,\gp {}_1],\,\G {}_2=(\gp {}_1,\gp {}_2],\ldots ,
\G {}_{s+1}=(\gp {}_s,\infty ).$$ We say that maximal ideals $\gp
$ and $\q $ from a linear orbit are {\em equivalent} ($\gp \sim \q
$) if they belong either to a non-degenerate orbit or to some $\G
{}_i $.

{\bf Weight $C$-modules.}  An $C$-module $V$ is {\em weight} if
${}_DV$ is semi-simple, i.e.
$$V=\oplus {}_{\gp \in {\rm Max}(D)}\,V_{\gp }$$
where $V_{\gp }=\{v\in V:\gp v=0\}=\{$the sum of simple
$D$-submodules which are isomorphic to ${}_D(D/ \gp ) \}$, the
{\em component} of $V$ of  {\em weight} $\gp $. The {\em support}
Supp$(V)$ of the weight module $V$ is the set of maximal ideals
$\gp $ such that $V_\gp \neq 0$.

Since
$$
XV_{\gp }\subseteq V_{\s (\gp )}\,\,{\rm and}\,\,YV_{\gp
}\subseteq V_{\s {}^{-1}(\gp )},
$$
each weight $C$-module $V$ decomposes into the direct sum of
$C$-submodules (the {\em orbit decomposition}) 
\begin{equation}\label{V=orb}
V=\bigoplus \{\,V_{\OO }\,|\,\OO\,{\rm is \,an \,orbit}\,\},
\end{equation}
where $V_{\OO }=\bigoplus \{\,V_{\gp }\,|\,\gp \in \OO \,\}.$
Obviously, for each maximal ideal $\gp $ of $D$ the module

$$C(\gp ):=C/C \gp \simeq C\otimes {}_D D/ \gp =\bigoplus {}_{i\in \mathbb{Z}}\,
v_i\otimes D/ \gp $$ is weight $({}_D(v_i\otimes D/ \gp)\simeq D/
\s {}^i(\gp ))$ with support $\OO (\gp ) .$

Denote by $\hC $(weight) the set of isoclasses of simple weight
$C$-modules. Each simple weight $C$-module and each simple
$D$-torsion $C$-module is a homomorphic image of $C(\gp )$ for
some $\gp \in {\rm Max}(D)$, and so by (\ref{V=orb}),
 \begin{equation} \label{Ator=weight}
   \hC (D-{\rm torsion})=\hC (D-{\rm weight}).
\end{equation}
We denote by $\hC ({\rm weight, linear})$ and  $\hC ({\rm weight,
cyclic})$ the sets of isoclasses of simple weight $C$-modules with
support from a {\em linear} and a  {\em cyclic} orbit
respectively. Then the set of isoclasses of simple weight
$C$-modules is the following disjoint union:
$$ \hC ({\rm weight})=\hC ({\rm weight, linear})\cup \hC ({\rm weight,
cyclic}).$$ The ideal $(H)$ of the polynomial algebra $D:=\CK [H]$
generated by the element $H$ is a maximal ideal of $D$. By
(\ref{spnHl}), $\l_1^{q^n}\neq \l_1$ for all $n\geq 1$
($\l_1:=[1]^\frac{1}{q}=x-x^\frac{1}{q}$), therefore the orbit
$\OO (H)$ of the maximal ideal $(H)$ is an infinite orbit. This is
the only degenerate orbit, this orbit is infinite (i.e. linear),
and there are only two equivalence classes in $\OO (H)$:
$\G_-:=(-\infty , (H)]$ and $\G_+:= ((H), \infty )$.

Let $Lin/\sim $ denote the set of equivalence classes in $Lin$,
that is
$$ Lin/\sim=\{ \G_-, \G_+, \OO (\gp )\, | \, \; {\rm where}\;\;
\OO (\gp ) \;\; {\rm is\; a\; non-degenerate\;  linear\;
orbit}\}.$$
\begin{theorem}\label{wlins}
The map
$$ Lin /\sim \ra \hC(D-{\rm weight, linear}),\,\G \ra [L(\G )],$$
is a bijection with inverse  $[L]\ra {\rm Supp}\,  L$ (in
particular, ${\rm Supp}\,  L(\G )=\G )$ where
\begin{enumerate}
\item if $\G \in$ Linn  then $L(\G )=C/C\gp $, for any $\gp \in
\G$, \item if $\G =(-\infty ,(H) ]\in$ then $L_-:=L(\G_-)=C/(CH +C
X)$, \item if $\G =((H) ,\infty )\in$ then $L_+:=L(\G_+ )=C/(C\s
(H )+CY)$.
\end{enumerate}
$\dim_{\CK}(L(\G ))=| \G | =\infty$ for each $[L(\G )]\in
\hC(D-{\rm weight, linear})$.
\end{theorem}

{\it Proof}. This is a special case of \cite{BavVOYJA97},
Corollary 4.1. $\Box $

Let us give more detail about the simple modules just described.

The simple weight $C$-module $L_-=\bigoplus_{i\geq 0}\CK \bY^i$,
$\bY^i:=Y^i+CH+CX$,  where the action of the generators for the
algebra $C$ is given by the rule:
\begin{eqnarray*}
H\bY^0=0,  & H\bY^i=-(\l_1+\l_1^q+\cdots +\l_1^{q^{i-1}})\bY^i, \; i\geq 1,  \\
X\bY^0=0,  & X\bY^i=-(\l_1+\l_1^q+\cdots
+\l_1^{q^{i-1}})\bY^{i-1}, \; i\geq 1,\\
Y\bY^i=\bY^{i+1}, & i\geq 0.
\end{eqnarray*}

The simple weight $C$-module $L_+=\bigoplus_{i\geq 0}\CK \bX^i$,
$\bX^i:=X^i+C\s (H)+CY$,
 where the action of the generators for the algebra $C$  is given
by the rule:
\begin{eqnarray*}
H\bX^0=\l_1\bX^0,  & H\bX=2\l_1\bX, \;\; H\bX^i=(2\l_1+\l_1^\frac{1}{q}+\cdots +
\l_1^\frac{1}{q^{i-1}})\bX^i, \; i\geq 2,  \\
Y\bX^0=0,  & Y\bX =\l_1\bX^0, \;\;
Y\bX^i=(2\l_1+\l_1^\frac{1}{q}+\cdots + \l_1^\frac{1}{q^{i-2}})\bX^{i-1}, \; i\geq 2,\\
X\bX^i=\bX^{i+1}, & i\geq 0.
\end{eqnarray*}

The simple weight $C$-module $$L=L(\G )=C/C\gp = (\bigoplus_{i\geq
1}\CK \bY^i)\bigoplus \CK \overline{1}\bigoplus (\bigoplus_{i\geq
1}\CK \bX^i) $$ where $\overline{u}:=u+C\gp$,
$\overline{1}:=\bY^0=\bX^0$ and $\gp =(H-\l )\in \G \in Linn$,
\begin{eqnarray*}
X\bX^i=\bX^{i+1}, \;\;\;\;\; Y\bY^i=\bY^{i+1}, & i\geq 0,\\
X\bY^i=(\l
-\l_1-\l_1^q-\cdots -\l_1^{q^{i-1}})\bY^{i-1}, & i\geq 1, \\
Y\bX^i=(\l +\l_1+\l_1^\frac{1}{q}+\cdots
+\l_1^\frac{1}{q^{i-1}})\bX^{i-1}, & i\geq 2,\\
Y\bX =\l \overline{1}. &
\end{eqnarray*}

The formulas above can be simplified taking into account that, for
$n\geq 2$,
\begin{eqnarray*}
\l_1+\l_1^q+\cdots +\l_1^{q^n}=-x^\frac{1}{q}+x^{q^n}, \\
\l_1+\l_1^\frac{1}{q}+\cdots
+\l_1^\frac{1}{q^n}=-x^\frac{1}{q^{n+1}}+x.
\end{eqnarray*}

\begin{lemma}\label{rfn}
For $\l \in \CK$ and a natural number $n$, $\s^n (D(H-\l ))=D(H-\l
)$ iff $\l \in -x^\frac{1}{q}+\mathbb{F}_{q^n}$.
\end{lemma}

{\it Proof}. If $n=1$ then $\s (D(H-\l ))=D(H-\l )$ iff $f_1=0$
where $f_1:= \l^q-\l +\l_1=\l^q-\l -x^\frac{1}{q}+x$. The
polynomial $f_1$ of degree $q$ (in $\l $) has $q$ {\em distinct}
roots since its derivative $\frac{df_1}{d\l}=-1\neq 0$. It follows
from the equality $ f_1(\l +\mu )=f_1(\l )+\mu^q-\mu$ that if $\l
$ is a root of the polynomial $f_1$ then so is $\l +\mu $ for each
$\mu\in \Fq$. Clearly, $-x^\frac{1}{q}$ is a root of the
polynomial $f_1$. Therefore, $-x^\frac{1}{q}+\mathbb{F}_{q}$ are
the roots of $f_1$.

Similarly, if $n\geq 2$ then, by (\ref{spnHl}), $\s^n (D(H-\l
))=D(H-\l )$ iff $f_n=0$ where
$$ f_n:= \l^{q^n}-\l +\l_1+\l_1^q+\cdots +\l^{q^{n-1}}=  \l^{q^n}-\l
-x^\frac{1}{q}+x^{q^{n-1}}.$$ The polynomial $f_n$ of degree $q^n$
(in $\l $) has $q^n$ {\em distinct} roots since its derivative
$\frac{df_n}{d\l}=-1\neq 0$. It follows from the equality $ f_n(\l
+\mu )=f_n(\l )+\mu^{q^n}-\mu$ that if $\l $ is a root of the
polynomial $f_n$ then so is $\l +\mu $ for each $\mu\in
\mathbb{F}_{q^n}$. Clearly, $-x^\frac{1}{q}$ is a root of the
polynomial $f_n$. Therefore, $-x^\frac{1}{q}+\mathbb{F}_{q^n}$ are
the roots of $f_n$. $\Box $

Recall that the {\em M\"{o}bious function} $\mu : \mathbb{N}\ra \{
0, \pm 1\}$ is given by the rule: $\mu (1)=1$, $\mu (p_1, \ldots
p_r)=(-1)^r$ where $p_1, \ldots , p_r$ are distinct primes, and
$\mu (n)=0$, otherwise. Given a function $f$ on $\mathbb{N}$ and a
second function $g$ defined by the formula $g(n)=\sum_{d|n}f(d)$.
Then (it is well-known) 
\begin{equation}\label{finvf}
f(n)=\sum_{d|n}\mu (d)g(\frac{n}{d}).
\end{equation}
The {\em Euler function} $\v$ on $\mathbb{N}$ is defined by $\v
(1)=1$ and, for $n>1$, $\v (n)=$ the number of natural numbers $m$
that are co-prime to $n$ and $1\leq m<n$.
\begin{theorem}\label{crfn}
(Classification of finite orbits) Let $\OO =\OO_\l $ $(\l \in \CK
)$ be the orbit of the maximal ideal $D(H-\l )$ of the polynomial
algebra $D:=\CK [H]$ under the action of the cyclic group $\langle
\s \rangle$. Then
\begin{enumerate}
\item the orbit $\OO_\l$ contains a single element iff $\l \in
-x^\frac{1}{q}+\mathbb{F}_{q}$. So, there are exactly $q$ distinct
maximal $\s$-invariant ideals of the algebra $D$,  \item the orbit
$\OO_\l$ contains exactly $n\geq 2$  elements iff $\l \in
-x^\frac{1}{q}+(\mathbb{F}_{q^n}\backslash \cup_{m|n, m\neq
n}\mathbb{F}_{q^m})$. So, there are exactly
$o_n:=n^{-1}\sum_{d|n}\mu (d)q^\frac{n}{d}=n^{-1}\v (q^n-1)$
distinct orbits that contain exactly $n\geq 2$ elements; and
$o_n\geq n^{-1}q^n(1-\frac{1}{q-1})>0$.
\end{enumerate}
\end{theorem}

{\it Proof}. $1$. This evident (see Lemma \ref{rfn}).

$2$. By Lemma \ref{rfn}, the orbit $\OO_\l $ contains exactly
$n\geq 2$ elements iff $\l \in
-x^\frac{1}{q}+(\mathbb{F}_{q^n}\backslash \cup_{m|n, m\neq
n}\mathbb{F}_{q^m})$. Let $o_n$ be the number of distinct orbits
that contain exactly $n$ elements, then $f(n)=no_n$ is the number
of all maximal ideals of $D$ that lie on all the orbits that
contain exactly $n$ elements. By Lemma \ref{rfn}, the function
$g(n):= \sum_{d|n}f(d)$ is equal to
$|-x^\frac{1}{q}+\mathbb{F}_{q^n}|=|\mathbb{F}_{q^n}|=q^n$ where
in this proof $|S|$ means the number of elements in a set $S$.
Therefore, by (\ref{finvf}), we have  $f(n)=n^{-1}\sum_{d|n}\mu
(d)q^\frac{n}{d}$. Clearly,
\begin{eqnarray*}
no_n &\geq &  |-x^\frac{1}{q}+\mathbb{F}_{q^n}|-|-x^\frac{1}{q}+
\mathbb{F}_{q^{n-1}}|-\cdots - |-x^\frac{1}{q}+\mathbb{F}_{q}|\\
 &=& q^n-q^{n-1}-\cdots -q= q^n-q\frac{q^{n-1}-1}{q-1}\\
 &> & q^n-\frac{q^n}{q-1}=q^n(1-\frac{1}{q-1})\geq 0.
\end{eqnarray*}
Clearly, $o_n=n^{-1}|\mathbb{F}_{q^n}\backslash \cup_{m|n, m\neq
n}\mathbb{F}_{q^m}|=n^{-1}\v (q^n-1)$. $\Box $

\begin{corollary}\label{1crfn}
 $\hC (D- {\rm torsionfree})=\emptyset$.
\end{corollary}

{\it Proof}. By \cite{BavVOYJA97}, Theorem 5.14, $\hC (D- {\rm
torsionfree})=\emptyset$ iff there are {\em infinitely} many {\em
cyclic}  orbits, and the result follows from Theorem
\ref{crfn}.(2). $\Box $

{\bf $\hC ({\rm weight, cyclic})$, simple finite dimensional (over
$\CK$) $C$-modules}. For a natural number $n$, set
$$C_{[n]} =\bigoplus_{i\in \mathbb{Z}} C_{in},$$
the {\em Veronese} subring of $C$. The ring $C$ considered as a
left (or right ) $C_{[n]} $-module is $(\Zn
=\mathbb{Z}/n\mathbb{Z})$-graded,
$$C={\bigoplus }_{\i \in \Zn }\, C_\i ,\;
{\rm where}\; C_\i ={\bigoplus }_{j\in \mathbb{Z} }\,C_{i+nj},\;\i
=i+n\mathbb{Z} .$$ It means that $C_{[n]} C_\i \subseteq C_\i $
for all $\i \in \Zn $. Moreover, it is a
 $\Zn $-graded ring (i.e. $C_\i C_\j \subseteq C_{\i +\j }$ for all
$\i ,\j \in \Zn )$.

Let $\OO $ be a {\em cyclic} orbit which contains $n=|\OO |$
elements. For a maximal ideal $\gp \in \OO $ $(\sn (\gp )=\gp )$,
let us consider the factor ring
$$C_{[n],\gp } :=C_{[n]} /(\gp ) $$
 of $C_{[n]} $ modulo  the ideal $(\gp )=\bigoplus_{i\in
\mathbb{Z}}(\gp v_{in} = v_{in} \gp)$ of $C_{[n]}$ generated by
$\gp $. Since every cyclic orbit is a non-degenerate one, the
algebra $C_{[n], \gp }=\CK [X^n, X^{-n}; \sn ]$ is the skew
Laurent extension with coefficients from the field $\CK$. Since
$\CK$ is a field and $\sn$ is an automorphism of the field $\CK$,
the ring $C_{[n],\gp }$ is a left principal ideal domain and a
right principal ideal domain. This means that every left (and
right) ideal of the ring $C_{[n],\gp }$  is generated by a single
element. The reason for this is that in the ring $C_{[n],\gp }$
the left (and right) division algorithm with remainder holds. So,
any left simple $C_{[n],\gp }$ -module has the form $C_{[n],\gp }/
C_{[n],\gp }b$ where $b$ is an irreducible element of $C_{[n],\gp
}$.

Let $M$ be a simple $C$-module with support from a cyclic orbit,
say $\OO$, that contains $n=|\OO |$ elements. It follows from the
weight decomposition
$$ M=\bigoplus_{\gp\in \OO }M_\gp$$
that $M_\gp$ is a simple $C_{[n],\gp }$-module for each $\gp \in
\Supp (M)$, and
$$ C_\i M_\gp\subseteq M_{\si (\gp )}\;\; {\rm for\; all} \;\; \i
\in \Zn .$$ It follows that $M\simeq C\t_{C_{[n],\gp }}M_\gp$ for
each $\gp \in \Supp (M)$, and so we have the following theorem.

\begin{theorem}\label{Cscyc}
\begin{enumerate}
\item $$\hC ({\rm weight, cyclic})=\cup_{n\geq 1}\cup_{\OO \in
Cyc_n} \hC (\OO )$$ where $\hC (\OO )$ is the set of isoclasses of
simple weight $C$-modules with support from the orbit $\OO $, and
$Cyc_n$ is the set of all the cyclic orbits that contain exactly
$n$ elements. \item For each orbit $\OO \in Cyc_n$ and a fixed
element $\gp \in \OO $, the map
$$ \hC_{[n],\gp} \ra \hC (\OO ), \;\; [N]\mapsto [C\t_{C_{[n], \gp}}N], $$ is a
bijection where $\hC_{[n],\gp}$ is the set of isoclasses of simple
$\hC_{[n],\gp}$-modules.
\end{enumerate}
\end{theorem}

Theorems \ref{wlins} and \ref{Cscyc} classify the simple
$C$-modules.

{\bf Finite dimensionality of kernels and cokernels}. In
\cite{MRJA73}, McConnel and Robson proved that for a simple module
$M$ over the first Weyl algebra $A_1(F)$ over a field $F$ of
characteristic zero and for any non-scalar element $u$ of
$A_1(F)$,  $\dim_F(\ker \, u_M)<\infty$ and  $\dim_F(\coker \,
u_M)<\infty$ where $u_M:M\ra M$, $m\mapsto um$. This result is
also true for some GWAs see \cite{BavGWA92}, and for the Carlitz
algebras.

\begin{theorem}\label{fdkck}
Let $M$ be a simple $C$-module and $u\in C\backslash \CK$. Then
$\dim_{\CK}(\ker \, u_M)<\infty$ and  $\dim_{\CK}(\coker \,
u_M)<\infty$.
\end{theorem}

{\it Proof}. The result is obvious if the module $M$ is finite
dimensional over $\CK$. So, let us assume that the module $M$ is
infinite dimensional over $\CK$, that is $M$ is from Theorem
\ref{wlins}. The $C$-module $M$ is a $\mathbb{Z}$-graded module,
each graded component is a $1$-dimensional vector space over the
field $\CK$ (see the description of modules after Theorem
\ref{wlins}).

If $u\in D$ then the kernel and cokernel of the map $u_M$
coincide,  and their common dimension over $\CK$ does not exceed
the number of distinct roots of the polynomial $u\in D:=\CK [H]$.

If $u=\alpha X^i$ or $u=\alpha Y^i$ for some $i\geq 1$ and $0\neq
\alpha \in D$ then the dimension of the kernel and cokernel of the
map $u_M$  over $\CK$ does not exceed
 $i +$ the number of distinct roots of the polynomial $\alpha$.

Finally, if $u=\alpha v_n+\beta v_{n-1}+\cdots +\g v_m$, $n>m$,
$\alpha $ and $\g $ are nonzero polynomials of $D$. Since
$$ \ker_\CK \, u_M\subseteq \ker_\CK \, (\alpha v_n)_M +\ker_\CK\,  (\g v_m)_M,$$
we have $\dim_{\CK} (\ker \, u_M)<\infty$.

Note that $M=\bigoplus_{i\in \mathbb{Z}}M_i$ is a
$\mathbb{Z}$-graded $C$-module where $M_i$ is a weight component
of dimension $1$ over $\CK$. Let $k=|n|+|m|$ and $M=\cup_{j\geq
0}M^j$ where $M^j:= \bigoplus_{-kj\leq i\leq kj}M_i$. Since
$\dim_{\CK} (\ker \, u_M)<\infty$, $uM^j\subseteq M^{j+1}$ for all
$j\geq 0$, and $\dim_{\CK}(M^{j+1})-\dim_{\CK}(M^{j})=2k=const$,
we have, for all $j\geq 0$,
\begin{eqnarray*}
\dim_{\CK}(M^{j+1}/uM^j) & =& \dim_{\CK}(M^{j+1})-\dim_{\CK}(uM^{j})\\
 &\leq & \dim_{\CK}(M^{j+1})-(\dim_{\CK}(M^{j})-\dim_{\CK}(\ker\,
 u_M))\\
 &=&\dim_{\CK}(M^{j+1})-\dim_{\CK}(M^{j})+\dim_{\CK}(\ker\,
 u_M)\\
 &=&c=const<\infty.
\end{eqnarray*}
Therefore, $\dim_{\CK}(\coker\,
 u_M)<c<\infty$. $\Box $


\section{The prime spectrum  and the ideal structure of the Carlitz
algebras}\label{prCarl}

In this section, the following theorem will be proved that
completely describes the ideal structure of the Carlitz algebra.

Recall that $Cyc_n$ denote the set of all the finite orbits that
contain exactly $n$ elements. By Theorem \ref{crfn}.(2),
$0<|Cyc_n|<\infty $. For each orbit $\OO \in Cyc_n$ and any $\l $
such that $D(H-\l )\in \OO$, let 
\begin{equation}\label{aOO}
\alpha_\OO :=\prod_{i=0}^{n-1}\si (H-\l ).
\end{equation}
Clearly, $\alpha_\OO$ is a {\em central} element of the algebra
$C$ (since $\sn (H-\l )=H-f_n(\l )-\l =H-\l $, see the proof of
Lemma \ref{rfn}), $\alpha_\OO $ is a {\em monic} polynomial of
$\CK [H]$ which is {\em does not depend} on the choice of $\l $.
The polynomial $\alpha_\OO$ is {\em the only monic} polynomial
which is a generator for the ideal $\prod_{\gp\in \OO}\gp$ of the
algebra $D$.

\begin{theorem}\label{mCarid}
\begin{enumerate}
\item Every nonzero prime ideal of the Carlitz algebra $C$ is a
maximal ideal. \item Maximal ideals of $C$ commute, $\gm\gn
=\gn\gm$. \item Each nonzero ideal $I$ of $C$ is a unique finite
product of maximal ideals of $C$: $I=\prod_{\gm \in {\rm Max}
(C)}\gm^{n(\gm )}$ for some $n(\gm)\geq 0$ all but finitely many
$n(\gm )=0$ and  if $ I=\prod_{\gm \in {\rm Max} (C)}\gm^{n(\gm
)}=\prod_{\gm \in {\rm Max} (C)}\gm^{l(\gm )}$ then $n(\gm )=l(\gm
)$ for all $\gm$. So, all ideals commute. \item The map $Cyc\ra
{\rm Max}(C)$, $\OO \mapsto \gm_\OO:=C\prod_{\gp \in \OO}\gp
=C\alpha_\OO$, is a bijection with inverse $\gm\mapsto \Supp
(C/\gm )$. \item For each $\OO \in Cyc$, the factor algebra
$C/\gm_\OO\simeq M_n(\CK [t,t^{-1};\sn ])$ is the $n\times n$
 matrix algebra with coefficients from the skew Laurent extension
 $\CK [t,t^{-1};\sn ]$ where $n=|\OO |$, the number of
 element in the orbit $\OO $. The centre $Z(C/\gm_\OO
 )\simeq \mathbb{F}_{q^n}$ and $\Kdim (C/\gm_\OO )=\gldim (C/\gm_\OO
 )=1$.
 \item The factor algebra $C/\gm_\OO$ is a domain iff $\OO =\{ \gp
 \}$ where $\gp$ is a $\s$-invariant maximal ideal of the algebra
 $D$. There are exactly $q$ such ideals (Theorem \ref{crfn}).
\item The factor algebras $C/\gm_\OO $ and $C/\gm_{\OO'}$ are
isomorphic over $\mathbb{F}_p$
 iff $|\OO |=| \OO'|$.
\end{enumerate}
\end{theorem}

In order to prove Theorem \ref{mCarid}, we first establish some
preliminary results that are interesting on their own.

The next result describes the centre of the algebra $C$ and its
localization $B$.
\begin{lemma}\label{Rs2}
\begin{enumerate}
\item For each natural number $n$, the skew Laurent extension $\CK
[t,t^{-1};\sn ]$ is a simple $\mathbb{F}_q$-algebra with centre
$\mathbb{F}_{q^n}$.
 \item The algebra $B=\S1 C=\CK (H)[X,X^{-1}; \s ]$ is a simple algebra with centre
 $Z(B)=\{ \mathbb{F}_{q}\prod_{\OO \in Cyc}\alpha_\OO^{n(\OO
 )}\, | \, n(\OO )\in \mathbb{Z}$ and all but finitely many  $n(\OO
 )=0\}$ that contains countably many elements.
 \item The centre $Z(C)=\{ \mathbb{F}_{q}\prod_{\OO \in Cyc}\alpha_\OO^{n(\OO
 )}\, | \, n(\OO )\geq 0$ and  all but finitely many  $n(\OO
 )=0\}$ is a unique factorization domain that contains countably many elements.
\end{enumerate}
\end{lemma}

{\it Proof}. $1$. For each natural $i$, $\CK^{\s^{ni}}:=\{ \l \in
\CK \, | \, \l =\s^{ni}(\l ) l\l^{q^{ni}}\} =\mathbb{F}_{q^{ni}}$.
We claim that the centre $Z$ of the $\mathbb{F}_p$-algebra $R:=\CK
[t,t^{-1}; \sn ]$ belongs to the field $\CK$ since otherwise we
would have a nonzero central element of the form $z=\l t^m+\cdots$
for some $0\neq m\in \mathbb{Z}$ and $0\neq \l \in \CK$ where the
three dots mean elements of strictly higher or strictly lower
degree in $t$. For any $\mu \in \CK \backslash
\mathbb{F}_{q^{n|m|}}$, $\mu z-z\mu =(\mu -\s^{nm}(\mu ))\l
t^m+\cdots \neq 0$, a contradiction. An element $\l \in \CK$
commutes with $t$ iff $\l \in \mathbb{F}_{q^n}$. Therefore,
$Z=\mathbb{F}_{q^n}$.

By \cite{MR}, 1.8.5, $R$ is a simple ring.

$2$. Similarly, by \cite{MR}, 1.8.5, $B$ is a simple ring. By
exactly the same reason as in the first case, $Z(B)\subseteq \CK
(H)$. Clearly, $\mathbb{F}_q\prod_{\OO \in Cyc}\alpha_{\OO}^{n(\OO
)}\subseteq Z(B)$. We have to prove that any nonzero element $z\in
Z(B)$ can be written  in this form. The rational function $z$ is
equal to $\g \frac{f}{g}$ where $f,g\in \CK [H]$ are co-prime
monic polynomials and $\g\in \CK$.

If $z=\g$ then $\g =X\g X^{-1}=\s (\g )$ , and so $\g \in
\mathbb{F}_q$, and we are done.

So, suppose that $z\neq \g $. Since $z=X^nzX^{-n}=\sn (z)$ for all
$n\in \mathbb{Z}$, and since $f$ and $g$ are co-prime we see that
$f$ and $g$ are equal to finite products of the form $\prod
\alpha_\OO^{n( \OO )}$ with $n(\OO )\geq 0$. Then $\g \in Z(B)$,
and so $\g \in \Fq$.

$3$. Clearly, $Z(C)\subseteq Z(B)$, and the result follows from
statement 2.  $\Box $

Recall that $Cyc_n=\{ \OO \in Cyc\, | \, | \OO | =n\}$.

\begin{lemma}\label{faCa}
For each $\OO \in Cyc_n$,
\begin{enumerate}
\item the ideal $C\alpha_\OO =\alpha_\OO C$ of $C$ generated by
the central element $\alpha_\OO$ of $C$ is a maximal ideal, and
 \item the factor ring $C/C\alpha_\OO \simeq M_n(\CK [t,t^{-1}; \sn ])$,
 the matrix algebra with entries from the skew Laurent extension
 $\CK [t,t^{-1}; \sn ]$.
  \item The centre $Z(C/C\alpha_\OO )\simeq F_{q^n}$.
\end{enumerate}
\end{lemma}

{\it Proof}. There are obvious $C$-module isomorphisms (where $\gp
\in \OO$):
\begin{eqnarray*}
C/C\alpha_\OO & \simeq & C\t_D D/D\alpha_\OO \simeq C\t_D
(\oplus_{i=1}^nD/\si (\gp ))\simeq \oplus_{i=1}^nC\t_DD/\si (\gp )\\
 &\simeq & \oplus_{i=1}^n C\t_DD/\gp \simeq (C\t_DD/\gp )^n.
\end{eqnarray*}
Now, there are obvious $\Fq$-algebra isomorphisms:
\begin{eqnarray*}
C/C\alpha_\OO & \simeq & {\rm End}_{C/C\alpha_\OO }(C/C\alpha_\OO
)\simeq {\rm End}_{C}((C\t_DD/\gp )^n)\simeq M_n({\rm
End}_{C}(C\t_DD/\gp ))\\
&\simeq & M_n({\rm End}_{C_{[n]}}(C_{[n], \gp}))\simeq M_n({\rm
End}_{C_{[n], \gp }}(C_{[n], \gp}))\simeq M_n(C_{[n], \gp})\simeq
M_n(\CK [t,t^{-1}; \sn ]).
\end{eqnarray*}
This proves the second statement. By Lemma \ref{Rs2},
$C/C\alpha_\OO$ is a simple algebra, hence $C\alpha_\OO$ is a
maximal ideal. The third statement follows from Lemma
\ref{Rs2}.(1),
$$ Z(C/C\alpha_\OO )\simeq Z(M_n(\CK [t,t^{-1};\sn ]))\simeq
Z(\CK [t,t^{-1};\sn ])\simeq \mathbb{F}_{q^n}. \;\;\; \Box $$

The next corollary describes the annihilators of simple
$C$-modules.

\begin{corollary}\label{ansimC}
\begin{enumerate}
\item ${\rm ann}_C(L)=0$ for any $[L]\in \hC ({\rm weight,
linear})$.\item ${\rm ann}_C(L)=C\alpha_{\Supp (L)}$ for any
$[L]\in \hC ({\rm weight, cyclic})$.
\end{enumerate}
\end{corollary}

{\it Proof}. $1$. This follows immediately from Theorem
\ref{wlins} and the description of the modules that follows
Theorem \ref{wlins}.

$2$. Let $\OO =\Supp (L)$. Then clearly, $C\alpha_\OO L=0$, and by
maximality of the ideal $C\alpha_\OO$ (Lemma \ref{faCa}.(1)) we
must have ${\rm ann}_C(L)=C\alpha_\OO$.  $\Box $

{\bf Proof of Theorem \ref{mCarid}}. $1$. Let $I$ be a proper
ideal of the algebra $C$. Since the algebra $B=\S1 C$ is a simple
algebra (Lemma \ref{Rs2}.(2)), we must have $\S1 I=B$, and so the
$C$-module $\bC :=C/I$ is a torsion one $(\S1 \bC = \S1
(C/I)\simeq \S1 C/\S1 I=0$). Its orbit decomposition
$$ \bC =\bigoplus_{\OO }\bC^\OO , \;\; \bC^\OO :=\ann_{\bC }(\alpha_\OO^i),$$
contains only finitely many summands (as $\bC$ is a Noetherian
module) and all the $\OO $ must be cyclic (since otherwise, there
exists $\OO \in Lin$, then $\bC_\OO$ contains a simple submodule,
say $L$, with support from an equivalence class of $\OO$, then
$0\neq I= {\rm ann}_C(\bC )\subseteq {\rm ann}_C(L)=0$, by
Corollary \ref{ansimC}, a contradiction). Suppose that $\OO_1,
\ldots , \OO_s$ are the cyclic orbits involved in the orbit
decomposition above and $\alpha_1:=\alpha_{\OO_1}, \ldots ,
\alpha_s:=\alpha_{\OO_s}$ are the corresponding central
polynomials. There exists a natural number $n$ such that
$\alpha_i^n\bC^{\OO_i}=0$ for all $i$. Therefore,
$C\prod_{i=1}^s\alpha_i^n\subseteq I$. If $I$ is a prime ideal
then $C\alpha_i\subseteq I$ for some $i$, and then $I=C\alpha_i$
since $C\alpha_i$ is a maximal ideal (Lemma \ref{faCa}.(1)).

$2$. We have just proved that any maximal ideal $\gm$ of the ring
$C$ is equal to the ideal generated by a central element
$\alpha_\OO$. Now, it is clear that maximal ideals commute.

$4$. Clearly, the map $Cyc\ra \Max (C)$, $ \OO \mapsto \gm_\OO
:=C\alpha_\OO$ is a bijection with inverse $\gm \mapsto
\Supp(C/\gm )$.

$5$ and $6$. Lemma \ref{faCa}, the (left and right) global
dimension $\gldim (C/\gm_\OO )=1$ (by \cite{MR}, 7.9.18) and the
(left and right) Krull dimension $\Kdim (C/\gm_\OO )=1$ (by
\cite{MR}, 6.5.4).

$7$. If $C/\gm_\OO \simeq C/\gm_{\OO'}$ then, by Lemma \ref{faCa},
(where $n=|\OO |$ and $n'=| \OO'|$)
$$ \mathbb{F}_{q^n}\simeq Z(C/\gm_\OO ) \simeq
Z(C/\gm_{\OO'})\simeq \mathbb{F}_{q^{n'}}, $$
 and so $n=n'$.

 If $n=n'$ then by Lemma \ref{faCa},
 $$ C/\gm_\OO\simeq M_n(\CK [t,t^{-1}; \sn ])\simeq
 C/\gm_{\OO'}.$$

$3$. Let $I$ be a proper ideal of the ring $C$, we keep the
notation of the proof of statement 1. The component $\bC^\OO$ is
equal to the union $\cup_{i\geq 1}\ann_{\bC}(\alpha^i_\OO)$, and
so $\bC^\OO$ is a two-sided ideal of the algebra $C$. Since
$\bC^\OO$ is a Noetherian $C$-module the chain
$$ \ann_{\bC}(\alpha_\OO) \subseteq \ann_{\bC}(\alpha_\OO^2)
\subseteq\cdots $$ must terminate, say on $n_\OO$ step, that is
$\bC^\OO =\ann_{\bC}(\alpha_\OO^{n_\OO})$. Since the polynomial
$\alpha_\OO\in D$ has only simple roots, it follows that for the
$C$-bimodule $C/C\alpha_\OO^i$ is a {\em uniserial} bimodule of
finite length, and the bimodule structure is given by the
descending chain of $C$-bimodules:
$$ C\supset C\alpha_\OO \supset C\alpha_\OO^2 \supset \cdots \supset
C\alpha_\OO^{i-1} \supset C\alpha_\OO^i$$
 and each subfactor $C\alpha_\OO^j/C\alpha_\OO^{j+1}\simeq
 C/C\alpha_\OO\simeq \CK [t,t^{-1};\sn ]$ $(n=|\OO |)$ is a simple
 algebra $=$ a simple $C$-bimodule (Lemma \ref{faCa}). Now, there
 exist unique numbers $n_i$ such that
 $$ \bC
 =\bigoplus_{i=1}^s\bC^{\OO_i}=
 \bigoplus_{i=1}^sC/C\alpha_i^{n_i}=C/C\prod_{i=1}^s\alpha_i^{n_i}.$$
Taking the annihilator we have
$$ I=\ann_C(\bC
)=\ann_C(C/C\prod_{i=1}^s\alpha_i^{n_i})=C\prod_{i=1}^s\alpha_i^{n_i}=
\prod_{i=1}^s\gm_{\OO_i}^{n_i}.$$ Note that $I\cap
D=D\prod_{i=1}^s\alpha_i^{n_i}$ and the uniqueness of the $n_i$
now is obvious.  Clearly, $I=C(I\cap D)=(I\cap D)C=C(I\cap Z(C))$.
$\Box$

In fact we have proved the following corollary.

\begin{corollary}\label{ct1}
\begin{enumerate}
\item If $I$ is an ideal of $C$ then $I=C(I\cap D)=(I\cap
D)C=C(I\cap Z(C))$. \item If $I$ and $J$ are  ideals of $C$ then
$I=J$ iff  $I\cap D=J\cap D$ iff $ I\cap Z(C)=J\cap Z(C)$.
\end{enumerate}
\end{corollary}


\section{The group of automorphisms $\Aut (C)$ and the isomorphism problem for  the Carlitz
algebras}\label{1kglC}

\begin{theorem}\label{caris}
Two distinct Carlitz rings are not isomorphic.
\end{theorem}

{\it Proof}. Given two distinct Carlitz rings $C_\nu $ and $C_\mu$
with $\nu <\mu$, and so $p^\nu <p^\mu$. Then, by Lemma
\ref{faCa}.(3) and Theorem \ref{mCarid}, the isomorphism invariant
numbers
\begin{eqnarray*}
p^\nu &=& \min \{ | Z(C_\nu /C_\nu \gm_\OO )|=p^{\nu |\OO |}\;\;\, | \, \OO\in Cyc\}, \\
p^\mu &=& \min \{ | Z(C_\mu /C_\mu \gm_{\OO'} )|=p^{\mu |\OO' |}\,
| \, \OO'\in Cyc\},
\end{eqnarray*}
are distinct. Therefore, the rings $C_\nu$ and $C_\mu$ cannot be
isomorphic.  $\Box $

For a domain $R$, let $R^*:=R\backslash \{ 0\}$. In particular,
$\CK^*$ is a multiplicative group.

\begin{theorem}\label{grrautC}
\begin{enumerate}
\item  The group of ring isomorphisms of the ring $C=C_\nu$, $\Aut
(C)=\{ \tau =\tau_{\alpha, \g , \d , \o }:X\mapsto \alpha X, \;
Y\mapsto \g \s^{-1}(\alpha^{-1})Y\; | \, $ where $\alpha \in
\CK^*$, $\g \in \mathbb{F}^*_{p^\nu}$, $\d\in \mathbb{F}_{p^\nu}$,
$\o \in \Aut (\CK)$ such that $\o (x)=\g x+\d\}$. \item The group
of $\CK$-isomorphisms of $C=C_\nu$, $\Aut_\CK (C)=\{ \tau
=\tau_{\alpha}:X\mapsto \alpha X, \; Y\mapsto
\s^{-1}(\alpha^{-1})Y\; | \, $ where $\alpha \in \CK^*\}\simeq
\CK^*$.
\end{enumerate}
\end{theorem}

{\it Proof}. $1$. Given $\tau \in \Aut (C)$. Then $\tau $ induces
the ring isomorphism of the centre of $C$. By Lemma \ref{Rs2}.(3)
and Theorem \ref{crfn}, we must have $\tau (H)=\g (H+\varepsilon
)$ for some $\g \in \Fq^*$ and $\varepsilon\in \Fq$. The algebra
$C$ is the GWA, and so $C=\oplus_{i\in \mathbb{Z}}C_i$ is a
$\mathbb{Z}$-graded algebra where $C_i=Dv_i$ is an eigenspace of
the inner derivation $\ad (H):C\ra C$, $c\mapsto [H, c]:=Hc-cH$,
that corresponds to the eigenvalue $e_i$ given by the rule
\begin{eqnarray*}
 \l_1+\l_1^q+\cdots +\l_1^{q^{i-1}}=-x^\frac{1}{q}+x^{q^{i-1}},& {\rm if} \;\; i\geq 1,  \\
  -\l_1-\l_1^\frac{1}{q}-\cdots -\l_1^\frac{1}{q^{|i|-1}}=-x+x^\frac{1}{q^{i}},&
   {\rm if} \;\; i\leq
  -1,\\
  0, & {\rm if} \;\; i=0,
\end{eqnarray*}
that is, $C_i=\{ c\in C\, | \, [H, c]=e_ic\}$. Note that all the
$e_i$ are distinct. Since $C_i=C_1^i$ and $C_{-i}=C_{-1}^i$ for
all $i\geq 1$. We must have either $\tau (C_{\pm 1})=C_{\pm 1}$
or, otherwise, $\tau (C_1)=C_{-1}$ and $\tau (c_{-1})=C_1$.

In the first case, $\tau (X)=\alpha X$ and $\tau (Y)=\beta Y$ for
some $\alpha , \beta \in D^*$. Applying $\tau $ to the identity
$YX=H$ we have the identity $\beta \s^{-1}(\alpha ) H=\g
(H+\varepsilon )$, therefore $\beta =\g \s^{-1}(\alpha^{-1})$, $
\alpha , \beta \in \CK^*$, and $\varepsilon =0$. Applying $\tau $
to the identity $YX-XY=\l_1$ and then dividing by $\g $ we get
$YX-XY=\g^{-1}\tau (\l_1)$, therefore, $\tau (\l_1)=\g \l_1$,
equivalently, $(\tau (x)-\g x)^q=\tau (x)-\g x$, and so $\tau
(x)=\g x+\d$ for some $\d \in \Fq$. The units of the ring $C$ form
the set $\CK^*$, therefore the $\tau$ induces an automorphism of
the field $\CK$, say $\o :=\tau |_\CK$. Thus $\tau =\tau_{\alpha,
\g , \d , \o }$. One can easily verify  that $\tau_{\alpha, \g ,
\d , \o } \in \Aut (C)$.

Let us show that the second case is impossible. In the second
case, $\tau (X)=\alpha Y$ and $\tau (Y)=\beta X$ for some $\alpha
, \beta \in D^*$. Applying $\tau $ to the identity $YX=H$ we have
$\g (H+\varepsilon )=\beta \s (\alpha )XY =\beta \s (\alpha )
(H-\l_1)$ which is impossible since $\varepsilon\in \Fq$ and
$\l_1\not\in \Fq$. So, the second case is vacuous.

$2$. Since $\o =\tau|_\CK =$ the identity map and $x\in \CK^*$, we
must have $\g =1$ and $\d =0$, and the result follows from the
first statement. $\Box $


\section{The Krull and the global dimensions of the Carlitz
algebras}\label{kglC}

{\bf The global dimension of the Carlitz algebra}.  The global
dimension, $\gldim$, means the left or right global dimension.

\begin{theorem}\label{gldimTHM}
(\cite{BavTHMA96}, Theorem 1.6) Let $A=D(\s , a)$ be a GWA where
$D$ is a commutative Dedekind domain, $Da=\gp_1^{n_1}\cdots
\gp_s^{n_s}$ be a product of distinct maximal ideals of $D$. Then
the global dimension of the algebra  $A$ is equal to

$$ \gldim (A)=\left\{%
\begin{array}{ll}
    \infty, & {\rm if} \; a=0\; {\rm or}\; \exists\;   n_i\geq 2, \\
    2, &{\rm if \; either } \; a\neq 0, \, n_1=\cdots =n_s=1, s\geq 1, \; {\rm or}\; a\;
    {\rm is \; invertible; \; and\; either }\\
    & \; Cyc\neq \emptyset\; {\rm or} \; \exists \, i\neq j \; {\rm s.t.} \;
    \s^k (\gp_i)=\gp_j\; {\rm for \; some}\; k\geq 1,
     \\
    1, & {\rm otherwise.} \\
\end{array}%
\right.
$$
\end{theorem}

\begin{corollary}\label{1gldimTHM}
The global dimension of the Carlitz algebra is $2$.
\end{corollary}

{\it Proof}. In the case of the Carlitz algebra $C$, we have
$D=\CK [H]$ and $a=H$, and $Cyc \neq \emptyset$. Therefore, by
Theorem \ref{gldimTHM}, $\gldim (C)=2$.  $\Box $

{\bf The Krull dimension of the Carlitz algebra}. The Krull
dimension, $\Kdim$, means the left or right Krull dimension.

Let $R$ be a  commutative Noetherian ring  and $\s \in \Aut (R)$.
A prime idea $\gp$ of $R$ is called $\s$-{\em semistable} if $\sn
(\gp )=\gp $ for some natural number $n\geq 1$. If there is no
such $n$ the ideal $\gp$ is called $\s$-{\em unstable}. ${\rm ht}
$ stands for the {\em height} of the prime ideal $\gp$.

\begin{theorem}\label{KdimGWA}
(\cite{BavVOYJA98}, Theorem 1.2) Let $R$ be a commutative
Noetherian ring with Krull dimension $\Kdim (R)<\infty$ and
 $A=R(\s , a)$ be a GWA. Then its Krull dimension $\Kdim (A)=\max
 \{ \Kdim (R), {\rm ht} \, \gp +1, {\rm ht}\, \gq +1\, |\, \gp$
 is a $\s$-unstable prime ideal of $R$ for which there exist
 infinitely many $i$ with $a\in \si (\gp )$; $\gq$ is a $\s$-semistable prime
  ideal of $R\}$.
\end{theorem}

\begin{corollary}\label{1KdimGWA}
The Krull dimension of the Carlitz algebra is $2$.
\end{corollary}

{\it Proof}. Note that $\Kdim (D)=1$ and so,  by Theorem
\ref{KdimGWA}, $\Kdim (C)\leq 2$. There exists a maximal ideal
$\gp$ of the algebra $D$ such that $\s (\gp )=\gp$. Then, by
Theorem \ref{KdimGWA},  we must have $\Kdim (C)=2$ since ${\rm
ht}\; \gp =1$.  $\Box $

Department of Pure Mathematics

University of  Sheffield

Hicks Building

Sheffield S3~7RH

UK

 email: v.bavula@sheffield.ac.uk

\end{document}